\documentclass[journal]{IEEEtran}
\usepackage{epsfig,graphicx,amssymb,amsmath}
\usepackage{bm}
\usepackage{amsthm}
\usepackage{amsmath,bm}
\usepackage{amsfonts}
\usepackage[noadjust]{cite}  
\usepackage{url}
\usepackage{array,tabularx}
\usepackage{indentfirst}
\usepackage{color}
\usepackage{stfloats}

\ifCLASSINFOpdf
\else
\fi
\begin{document}
%
\title{Distributed algorithms for the least square solution of linear equations}
%


\author{Yi~Huang
        and Ziyang~Meng
\thanks{This work has been supported in part by the National Natural Science Foundation of China under Grants
U19B2029, 61873140, and 61833009, and the Institute for Guo Qiang of Tsinghua University under Grant 2019GQG1023. (Corresponding author: Ziyang Meng.)

Yi Huang and Ziyang Meng are with the Department of Precision
Instrument, Tsinghua University, Beijing 100084, China (e-mail:
huangyi161531@126.com; ziyangmeng@mail.tsinghua.edu.cn)}}
\maketitle

\begin{abstract}
This paper proposes distributed algorithms for solving linear equations to seek a least square solution via multi-agent networks. We consider that each agent has only access to a small and imcomplete block of linear equations rather than the complete row or column in the existing literatures. Firstly, we focus on the case of a homogeneous partition of linear equations. A distributed algorithm is proposed via a single-layered grid network, in which each agent only needs to control three scalar states. Secondly, we consider the case of heterogeneous partitions of linear equations. Two distributed algorithms with doubled-layered network are developed, which allows each agent's states to have different dimensions and can be applied to heterogeneous agents with different storage and computation capability. Rigorous proofs show that the proposed distributed algorithms collaboratively obtain a least square solution with exponential convergence, and also own a solvability verification property, i.e., a criterion to verify whether the obtained solution is an exact solution. Finally, some simulation examples are provided to demonstrate the effectiveness of the proposed algorithms.
\end{abstract}

\begin{IEEEkeywords}
Distributed algorithms, least square solution, linear equations, multi-agent networks
\end{IEEEkeywords}

\newtheorem{Assumption}{Assumption}
 \newtheorem{Remark}{Remark}
 \newtheorem{Lemma}{Lemma}
 \newtheorem{Definition}{Definition}
  \newtheorem{Proposition}{Proposition}
 \newtheorem{Theorem}{Theorem}
 \newtheorem{Property}{Property}
 \newtheorem{Corollary}{Corollary}
 \newtheorem{Example}{Example}

\section{Introduction}
The problem of solving a linear algebraic equation has received significant consideration due to its broad application in practical engineering including wireless sensor network \cite{a1}, simultaneous localization and mapping \cite{a2}, electromagnetism computations \cite{a3}, PageRank problem of search engine \cite{a4}, and etc. Some classical algorithms including Jacobi method, Gauss-Seidel method, are first developed in a centralized manner \cite{1}. Nevertheless, due to the limited computation, storage and communication resources of the central node, these algorithms are not applicable to the case of large-scale linear equations with a large number of unknown variables. Compared with these centralized methods, distributed algorithms are more effective to solve the large-scale linear equations \cite{2}. By virtue of multi-agent networks, distributed algorithm decomposes the large-scale linear equations into many small-scale linear equations for each agent in the network. By exchanging the information with local neighbors, all the agents can solve the original linear equation in a collaborative way.

Based on consensus control and optimization method, many distributed algorithms have been proposed to solve the linear equations $Ax=b$. For instance, by considering each agent only has access to a few rows of $A$ and $b$, discrete-time distributed algorithms are proposed in \cite{3,4,5,6,6a} and continuous-time algorithms are developed in \cite{7,8,8a,8b}. However, most of these existing algorithms are proposed for solving the linear equations which have exact solutions. In the practical applications, many linear equations are over-determined where the number of linear equations is much larger than that of the unknown variables. In general, the over-determined linear equation has no exact solution and only seeks a least square solution. For the case of over-determined linear equations, above mentioned algorithms can not achieve a least square solution.

To seek a least square solution of linear equations, some distributed algorithms are proposed by transforming the least square problem into a distributed optimization problem. For instance, a continuous-time distributed algorithms is developed by using the classical Arrow-Hurwic-Uzawa flow to find a least square solution and its discrete version is given in \cite{9}. A consensus-based sub-gradient continuous-time algorithms is proposed to compute a least square solution in \cite{10}. By considering each agent knows a few columns of $A$ and $b$, a distributed algorithm by using saddle point dynamics is proposed such that a least square solution is obtained in \cite{10a}. In addition, several discrete-time distributed algorithms are developed in \cite{11a, 11,12} to obtain a least-square solution.

We note that most of above algorithms require that each agent knows at least one complete row information of the overall linear equations and the estimated state of each agent has same dimension as the unknown variables to be solved. For the case of the large-scale linear equations with large dimension of the unknown variables, above mentioned algorithms may fail by considering that each agent only has limited storage, computation and communication resources. Recently, in \cite{13}, a continuous-time distributed algorithm with scalar state is proposed for solving linear equation, in which each agent has only access to two scalar elements of $A$ and $b$, and the estimated states are also two scalars. Furthermore, some scalable distributed algorithms for solving linear equations are proposed by introducing double-layered multi-agent networks in \cite{14}. However, the proposed algorithms in \cite{13,14} are only effective for the linear equations that have exact solutions, and cannot be applied for over-determined linear equations to achieve a least square solution.

Motivated by above discussions, the objective of this paper is to develop distributed algorithms with small state dimension for each agent to find a least square solution of linear equations. Inspired by the work in \cite{14}, we consider that one homogeneous partition and two different heterogeneous partitions of linear equations, in which each agent has only access to a small and imcomplete block of the linear equations (can be as small as two scalar elements). Firstly, we develops a distributed algorithm via a single-layered grid network for homogeneous partition case, in which each agent only needs to control three scalar states. Secondly, two distributed algorithms via doubled-layered network for two heterogeneous partitions are proposed, in which each agent's states are allowed to be different dimensions and even scalars. Compared with the existing results, the main contributions of this paper are three-fold:

(1) Unlike the existing algorithms in \cite{7,8,8a,8b,9,10,10a,11a,11,12} that requires complete row or column information, we consider that each agent has only access to two scalar elements when the considered linear equation is partitioned as small as possible. In such a case, each agent only has three scalar states and therefore has low requirements of computation, storage and communication resources for solving the large-scale equations.

(2) In contrast with the distributed algorithms in \cite{13,14} that are only effective for the linear equations that have at least one solution, the proposed distributed algorithms can be applied for over-determined linear equations to achieve a least-square solution. In particular, the ideas of dynamic average consensus is introduced, and a crucial eigenvalue analysis of the general block matrix is provided to show the exponential convergence of the proposed algorithms.

(3) The proposed distributed algorithms do not require the assumption in \cite{9,12} that the least square solution of linear equation is unique. Most importantly, the proposed algorithms are capable of verifying whether the obtained least-square solution is an exact solution or not.

The rest of this paper is organized as follows. In Section II, we formulate the least square problem of linear equation. Several distributed algorithms under homogenous and heterogenous partitions are presented in Section III and IV, respectively. Finally, conclusions are drawn in Section V.

\section{Preliminaries and problem formulation}

Let $\mathbb{R}^{n\times n}$ and $\mathbb{R}^{n}$ be the sets of $n\times n$ real matrices and $n$ dimension real vectors. $\mathbb{C}$ denotes the complex space. Let $I_{n}$ be a $n\times n$ identity matrix and $1_{n}$ be a $n$-dimensional vector with all entries being 1. $(\cdot)^{T}$ represents the transpose in the real space, and $(\cdot)^{H}$ denotes the conjugate transpose in the complex space. $\Vert \cdot\Vert$ denotes the Euclidean norm of a vector or a matrix. Let $\mathcal{I}_{n}=\{1,2,\ldots, n\}$ be a set of natural numbers up to $n$. $\text{col}(x_{1}, \ldots, x_{n})$ denotes a column stack of vectors $x_{i}, i \in \mathcal{I}_{n}$ and $\text{diag}(x_{1},x_{2},\ldots, x_{n})$ denotes the diagonal matrix with diagonal entries being $x_{1}$ to $x_{n}$.  $A \otimes B$ refers to the Kronecker product of the matrices $A$ and $B$. $\text{ker}(M)$ and $\text{Image}(M)$ denote the kernel space and image space of matrix $M$, respectively.

Consider the following linear algebraic equation
\begin{align}\label{1}
Ax=b, x\in \mathbb{R}^{n},
\end{align}
where $A\in \mathbb{R}^{m\times n}$ and $b\in \mathbb{R}^{m}$. Note that if $b\in \text{Image}(A)$, linear equation \eqref{1} always has a solution. If $b \notin \text{Image}(A)$, equation \eqref{1} has no solution and its least square solution is defined by the following optimization problem
\begin{align}\label{2}
\text{min}_{x\in \mathbb{R}^{n}} \Vert Ax-b\Vert^{2}.
\end{align}

\begin{Lemma}[\cite{11}]
The variable $x^{*} \in \mathbb{R}^{n}$ is a least square solution of equation \eqref{1} if and only if $A^{T}(Ax^{*}-b)=0$.
\end{Lemma}

\section{Distributed algorithms design under Homogenous partition}
In this section, we consider a homogeneous partition of the linear equation \eqref{1}.  The homogeneous partition means that all the complete row blocks of matrix $A \in \mathbb{R}^{m\times n}$ have same column partitions or  all the complete column blocks also have same row partitions. To make the analysis more intuitive, we consider the homogeneous partition of \eqref{1} with $m$ row blocks and $n$ column blocks such that
\begin{align}\label{3}
A=\text{col}(A_{1},A_{2},\ldots, A_{m}), A_{i}=[A_{i1},A_{i2},\ldots, A_{in}],
\end{align}
with the scalar $A_{ij}, i\in \mathcal{I}_{m}, j\in \mathcal{I}_{n}$. Correspondingly, $b=\text{col}(b_{1},b_{2},\ldots, b_{m})$ and $b_{i}=\sum^{n}_{j=1}b_{ij}$ with the scalar $b_{ij}, i\in \mathcal{I}_{m}, j\in \mathcal{I}_{n}$. Note that $b_{i}$ is a partition of $b$ and $b_{ij}$ can be chosen arbitrarily by satisfying $b_{i}=\sum^{n}_{j=1} b_{ij}$.

\begin{figure}[!ht]
\centering
\includegraphics[width=0.25\textwidth, clip=true]{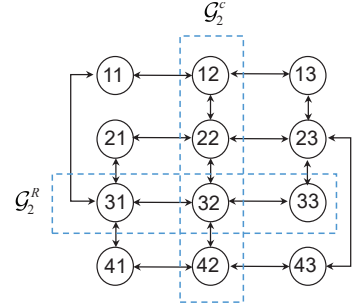}
\caption{Single-layered grid multi-agent network network}
\end{figure}

Inspired by the results given in \cite{14}, we introduce a single-layered grid multi-agent network graph $\mathcal{G}$ consisting of $mn$ agents. Each agent in the graph $\mathcal{G}$ is denoted as $ij, i\in  \mathcal{I}_{m}, j \in  \mathcal{I}_{n}$, and each agent $ij$ can communicate with some other agents called its neighbors. The communication relationship among the row agents $\mathcal{V}^{R}_{i}=\{i1, \ldots, in\}, \forall i\in \mathcal{I}_{m}$ is denoted by an undirected graph $\mathcal{G}^{R}_{i}=\{\mathcal{V}^{R}_{i}, \mathcal{E}^{R}_{i}\}$ with the edge set $\mathcal{E}^{R}_{i}\subseteq \mathcal{V}^{R}_{i}\times \mathcal{V}^{R}_{i}$. The row neighbors of agent $ij$ is denoted as $\mathcal{N}^{R}_{ij}=\{ik~|~ (ik,ij)  \in  \mathcal{E}^{R}_{i}\}, \forall i\in \mathcal{I}_{m}$. Similarly, the communication topology among the column agents $\mathcal{V}^{c}_{j}=\{1j,\ldots, mj\}, \forall j \in \mathcal{I}_{n}$ is denoted by the undirected graph $\mathcal{G}^{c}_{j}=(\mathcal{V}^{c}_{j},\mathcal{E}^{c}_{j})$ with the edge set $\mathcal{E}^{c}_{j}\subseteq \mathcal{V}^{c}_{j}\times \mathcal{V}^{c}_{j}$, and the column neighbors of agent $ij$ are denoted as $\mathcal{N}^{c}_{ij}=\{kj~|~ (kj,ij)  \in  \mathcal{E}^{c}_{j}\}, \forall j \in \mathcal{I}_{n}$. Fig. 1 shows an example of such a grid multi-agent network of homogeneous partition.

\begin{Assumption}
The graphs $\mathcal{G}^{R}_{i}$ and $\mathcal{G}^{c}_{j}$ for any $i\in \mathcal{I}_{m}, j\in \mathcal{I}_{n}$ are undirected and connected .
\end{Assumption}
Suppose that each agent $ij, i\in\mathcal{I}_{m}, j\in\mathcal{I}_{n}$ has only access to two scalars $A_{ij}, b_{ij}$. In this section, we target to propose a distributed algorithm such that the least square problem (2) is solved. In particular, the objective is to guarantee that each agent $ij$'s state $x_{ij}$ converges exponentially to a constant scalar $x^{*}_{ij}$ satisfying that

(i) All $x^{*}_{ij}$ in the same column achieve a consensus $x^{*}_{j}$, i.e.,
\begin{align}\label{4}
x^{*}_{1j}=x^{*}_{2j}=\ldots=x^{*}_{mj}=x^{*}_{j}, j \in \mathcal{I}_{n}.
\end{align}

(ii) The consensus values $x^{*}_{j}, j\in \bm n$ among all the columns satisfying that
\begin{align}\label{5}
\sum^{m}_{i=1}A_{ij}\sum^{n}_{j=1}(A_{ij}x^{*}_{j}-b_{ij})=0.
\end{align}
Note that $\sum^{n}_{j=1}A_{ij}x^{*}_{j}-b_{ij}=A_{i}x^{*}-b_{i}$ with $x^{*}=\text{col}(x^{*}_{1}$ $,\ldots, x^{*}_{n})$. From \eqref{5}, it follows that $\sum^{m}_{i=1}A_{ij}(A_{i}x^{*}-b_{i})=0, \forall j\in \bm n$, which implies that $A^{T}(Ax^{*}-b)=0$. Based on Lemma 1, we know that the constant value $x^{*}_{ij}, j\in \bm n$ collaboratively forms a least square solution of \eqref{1}.

\subsection{Distributed algorithm}
To ensure that all the agents' states $x_{ij}$ converge to $x^{*}_{ij}$ satisfying \eqref{4} and \eqref{5}, we introduce two extra variables $y_{ij},z_{ij} \in \mathbb{R}$ for each agent $ij$. Then, a distributed algorithm for each agent $ij, i\in \bm m, j\in \bm n$ is given as
\begin{align}\label{6}
\begin{cases}
\dot{x}_{ij}=-(A_{ij}y_{ij}-\sum_{kj \in \mathcal{N}^{c}_{ij}}(z_{ij}-z_{kj}))\\
~~~~~~~~-\sum_{kj \in \mathcal{N}^{c}_{ij}}(x_{ij}-x_{kj}),\\
\dot{\xi}_{ij}=-\sum_{ik \in \mathcal{N}^{R}_{ij}}(y_{ij}-y_{ik}),y_{ij}=\xi_{ij}+\\
~~~~~~~~~A_{ij}x_{ij}-b_{ij},\\
\dot{z}_{ij}=A_{ij}y_{ij}-\sum_{kj \in \mathcal{N}^{c}_{ij}}(z_{ij}-z_{kj}),
\end{cases}
\end{align}
where $\xi_{ij} \in \mathbb{R}$ is an intermediate variable and its initial value is selected to be zero. In fact, the second equation of \eqref{6} utilizes the ideas of dynamic average consensus approach in \cite{15,16} such that the average values $\frac{1}{n}\sum^{n}_{j=1}(A_{ij}x_{ij}-b_{ij})$.
can be estimated by $y_{ij}$ for all $i\in \bm n$. In addition, the term $-\sum_{kj \in \mathcal{N}^{c}_{ij}}(x_{ij}-x_{kj})$ in the dynamics of $x_{ij}$ is introduced to achieve the column consensus \eqref{4}, and the term $A_{ij}y_{ij}-\sum_{kj \in \mathcal{N}^{c}_{ij}}(z_{ij}-z_{kj})$ in the dynamics of $x_{ij}$ and $z_{ij}$ are used to guarantee that \eqref{5} is satisfied.

Define $\bar{x}_{i}=\text{col}(x_{i1},\ldots, x_{in}), \bar{\xi}_{i}=\text{col}(\xi_{i1},\ldots, \xi_{in}), \bar{y}_{i}=\text{col}(y_{i1}, \ldots, y_{in}), \bar{z}_{i}=\text{col}(z_{i1},\ldots, z_{in}) \in \mathbb{R}^{n}$. Let $\bar{A}_{i}=\text{diag}\{A_{i1},\ldots, A_{in}\} \in \mathbb{R}^{n\times n}$, and $L_{\mathcal{G}^{R}_{i}}\in \mathbb{R}^{n\times n}$ is the Laplacian matrix of the graph $\mathcal{G}^{R}_{i}$. Then, for each $i\in \bm m$, the algorithm \eqref{6} can be written in the following compact form
\begin{align}\label{7}
\begin{cases}\dot{\bar{x}}_{i}=-(\bar{A}_{i}\bar{y}_{i}-\bar{\psi}_{i})-\bar{\phi}_{i},\\
\dot{\bar{y}}_{i}=-(L_{\mathcal{G}^{R}_{i}}+\bar{A}_{i}\bar{A}^{T}_{i})\bar{y}_{i}+\bar{A}_{i}\bar{\psi}_{i}-\bar{A}_{i}\bar{\phi}_{i},\\
\dot{\bar{z}}_{i}=\bar{A}_{i}\bar{y}_{i}-\bar{\psi}_{i},
\end{cases}
\end{align}
where $\bar{\psi}_{i}=\text{col}(\psi_{i1},\ldots, \psi_{in})\in \mathbb{R}^{n}$ with $\psi_{ij}=\sum_{kj\in \mathcal{N}^{c}_{ij}}(z_{ij}$ $-z_{kj})$, $\bar{\phi}_{i}=\text{col}(\phi_{i1},\ldots, \phi_{in})\in \mathbb{R}^{n}$ with $\phi_{ij}=\sum_{kj\in \mathcal{N}^{c}_{ij}}(x_{ij}$ $-x_{kj}), i\in \bm m, j\in \bm n$. According to the definition of the Laplacian matrix $L_{\mathcal{G}^{c}_{j}}\in \mathbb{R}^{m\times m}$ of $\mathcal{G}^{c}_{j}$, one can derive that
\begin{align}\label{8}
\begin{bmatrix} \psi_{1j} \\ \psi_{2j} \\ \ldots \\ \psi_{mj} \end{bmatrix}=L_{\mathcal{G}^{c}_{j}}\begin{bmatrix} z_{1j}\\ z_{2j} \\ \ldots \\  z_{mj} \end{bmatrix}, \begin{bmatrix} \phi_{1j} \\ \phi_{2j} \\ \ldots \\ \phi_{mj} \end{bmatrix}=L_{\mathcal{G}^{c}_{j}}\begin{bmatrix} x_{1j}\\ x_{2j} \\ \ldots \\  x_{mj} \end{bmatrix}.
\end{align}
Define $\hat{x}=\text{col}(\bar{x}_{1},\ldots, \bar{x}_{m}), \hat{z}=\text{col}(\bar{z}_{1},\ldots, \bar{z}_{m}) \in \mathbb{R}^{mn},$ $ \hat{\psi}=\text{col}(\bar{\psi}_{1},\ldots, \bar{\psi}_{m})\in \mathbb{R}^{mn}$ and $\phi=\text{col}(\bar{\phi}_{1},\ldots, \bar{\phi}_{m})\in \mathbb{R}^{mn}$. Let $P\in \mathbb{R}^{mn \times mn}$ be a row permutation matrix with $P^{T}P=I_{mn\times mn}$ such that
\begin{align}\label{9}
P\hat{x}=\text{col}(x_{11}, x_{21},\ldots, x_{m1},\ldots, x_{1n},x_{2n},\ldots, x_{mn}).
\end{align}
It then from \eqref{8} and \eqref{9} that
\begin{align}\label{9a}
P\hat{\psi}=\hat{L}_{c}P\hat{z},  P\hat{\phi}=\hat{L}_{c}P\hat{x},
\end{align}
with $\hat{L}_{c}=\text{diag}(L_{\mathcal{G}^{c}_{1}},\ldots, L_{\mathcal{G}^{c}_{n}})\in \mathbb{R}^{mn\times mn}$. Let $\hat{A}=\text{diag}(\bar{A}_{1}$ $,\ldots, \bar{A}_{m}), \hat{y}=\text{col}(\bar{y}_{1},\ldots, \bar{y}_{m}), \hat{L}_{R}=\text{diag}(L_{\mathcal{G}^{R}_{1}},$ $\ldots, L_{\mathcal{G}^{R}_{m}})$ and $\tilde{L}_{c}=P^{T}\hat{L}_{c}P$. It then follows that \eqref{7} can be written as
\begin{align}\label{10}
\begin{cases}
\dot{\hat{x}}=-(\hat{A}\hat{y}-\tilde{L}_{c}\hat{z})-\tilde{L}_{c}\hat{x},\\
\dot{\hat{y}}=-(\hat{L}_{R}+\hat{A}\hat{A})\hat{y}+\hat{A}\tilde{L}_{c}(\hat{z}-\hat{x}),\\
\dot{\hat{z}}=\hat{A}\hat{y}-\tilde{L}_{c}\hat{z},
\end{cases}
\end{align}
which is equivalent as
\begin{align}\label{11}
\begin{bmatrix} \dot{\hat{x}} \\ \dot{\hat{y}} \\ \dot{\hat{z}}\end{bmatrix}=\begin{bmatrix} -\tilde{L}_{c}& -\hat{A}& \tilde{L}_{c}\\
-\hat{A}\tilde{L}_{c}& -(\hat{L}_{R}+\hat{A}\hat{A}) & \hat{A}\tilde{L}_{c} \\
0& \hat{A}& -\hat{A}\tilde{L}_{c} \end{bmatrix}\begin{bmatrix} \hat{x} \\ \hat{y} \\ \hat{z} \end{bmatrix}.
\end{align}
To analyze the convergence of system \eqref{11}, we provide the following Lemma to characterize the eigenvalues of the system matrix of \eqref{11}.
\begin{Lemma}
Let
\begin{align}\label{12}
M=\begin{bmatrix} -L_{1} & -Y & L_{1}\\ -Y^{T}L_{1} & -Y^{T}Y-L_{2} & Y^{T}L_{1}\\ 0 & Y & -L_{1} \end{bmatrix}
\end{align}
where $L_{1}, L_{2}, Y$ are real and positive semi-definite, $Y$ is a real matrix with appropriate dimension. Then, we can obtain that

(i) any non-zero eigenvalues of $M$ has negative real part.

(ii) zero eigenvalue must be non-defective, that is, its algebraic multiplicity is equal to its geometric multiplicity.
\end{Lemma}

The proof of Lemma 2 can be seen in Appendix A. Based on Lemma 2, the convergence result of the system \eqref{11} is given in the following theorem and its proof can be seen in Appendix B.
\begin{Theorem}
Suppose that Assumption 1 holds and the initial value of $\xi_{ij}, i\in \bm m, j\in \bm n$ selected to be zero. The proposed distributed algorithm \eqref{6} guarantees that all the agents' state $x_{ij}(t), i\in \bm m, j\in \bm n$ converge exponentially to a constant value $x^{*}_{ij}$ satisfying \eqref{4} and \eqref{5}, which collaboratively forms a least square solution $x^{*}$ of equation \eqref{1}. Moreover, $x^{*}$ is an exact solution of \eqref{1} if and only if $\lim_{t\to \infty} y_{ij}(t)=0$.
\end{Theorem}

\begin{Remark}
Note that the distributed algorithms in \cite{14} can be only applied to the determined linear equation with an exact solution and cannot be used to the over-determined linear equation. In contrast, we propose a new distributed algorithm to obtain a least square solution of determined linear equation. Compared with the results in \cite{14}, the main technical difficulties are given as follows: (1) in general, seeking a least square solution is more complex than an exact solution. To achieve a least square solution in a distributed manner, we utilize the dynamic average consensus approach to estimate the sum term $\frac{1}{n}\sum^{n}_{j=1}(A_{ij}x_{ij}-b_{ij})$ through the variable $y_{ij}$. (2) Since one extra variable is introduced, the system matrix of the proposed algorithm \eqref{6} is three-by-three block matrices rather than two-by-two block matrices in \cite{14}. This leads to that the eigenvalue analysis of the proposed algorithm \eqref{6} is more difficult. To deal with this issue, we provide a detailed eigenvalue analysis for more general block matrices shown in Lemma 2, which is crucial to show the exponential convergence of the proposed algorithm.
\end{Remark}

\begin{Remark}
For the problem of solving linear equation $Ax=b$, the distributed algorithms in \cite{7,8,8a,8b,9,10,10a,11a,11,12} that requires complete rows or columns for each agent are involved with vector operation rather than matrix operation such that the computation complexity for each agent can be reduced. In contrast with the existing algorithms given in \cite{7,8,8a,8b,9,10,10a,11a,11,12}, we propose scalable distributed algorithm \eqref{6} only involving with scalar operation for each agent rather than vector operation. Thus, the computation, storage and communication burden of each agent can be reduced effectively and can be applied to the large-scale linear equations.
\end{Remark}

\subsection{Simulation}
In this section, some simulation examples are provided to verify the effectiveness of the proposed algorithm \eqref{6} for solving $Ax=b$ under homogenous partitions in \eqref{3}, where $Ax=b$ is given as
\begin{align}\label{s1}
A=\begin{bmatrix}
    1 & 2 & 1 \\
    2 & -1 &  -1\\
    1 & -2 &  4\\
    2 & 2  & -2\end{bmatrix}, b=\begin{bmatrix} 3\\ 2\\ 1\\ 2 \end{bmatrix}.
\end{align}
Since $\text{rank}(A)=3<\text{rank}(A,b)=4$, it follows that above linear equation $Ax=b$ has no exact solution and its least square solution is given as $x^{*}=(A^{T}A)^{-1} A^{T}b=[1.1310, 0.4947, 0.2992]^{T}$. From \eqref{3}, one can derive that $A_{ij}$ is the $ij$th entry of $A$. In addition, we choose $b_{ij}=b_{i}$ if $j=i$ while $b_{ij}=0$ if $j\ne i$. We utilize the grid multi-agent network $\mathcal{G}$ in Fig. 1 to solve collaboratively the solution $x^{*}$.  Let $E_{1}(t)=\sum^{m}_{i=1}\sum^{n}_{j=1}\Vert x_{ij}(t)-x^{*}_{i}\Vert^{2}, x^{*}=\text{col}(x^{*}_{1},\ldots, x^{*}_{n})$ be the total estimated error, and $Y_{e}(t)=\sum^{m}_{i=1}\sum^{n}_{j=1}\Vert y_{ij}(t)\Vert^{2}$ is equal to $\Vert \hat{y}(t)\Vert^{2}$. Fig. 2 shows that $E_{1}(t)$ converges exponentially fast to zero. It is illustrated that all the states $x_{ij}(t)$ converge exponentially to the constant value $x^{*}_{ij}$, which consists of the solution $x^{*}=\text{col}(x^{*}_{i1}, \ldots, x^{*}_{in})$. Fig. 2 also shows that $Y_{e}(t)$ finally converges a non-zero value, which implies that $\lim_{t\to \infty} \hat{y}(t) \ne 0$. Based on the results of Theorem 1, we can verify that the obtained solution $x^{*}$ is not an exact solution of \eqref{s1}.

\begin{figure}[!ht]
\centering
\includegraphics[width=0.45\textwidth, clip=true]{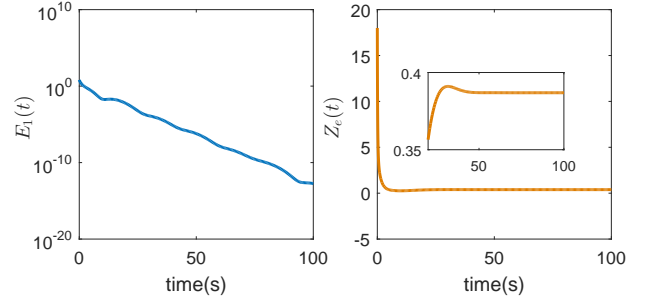}
\caption{Evolutions of $E_{1}(t)$ and $Y_{e}(t)$ under algorithm \eqref{6}}
\end{figure}

\section{Distributed algorithm design under heterogeneous partitions}
In this section, we consider heterogeneous partitions of \eqref{1}, which includes two cases: Case 1 is that each complete row block of $A$ has different column partitions, and Case 2 is that each complete column has different row partitions. One example of these two partition cases is shown in Fig. 3.
\begin{figure}[!ht]
\centering
\includegraphics[width=0.35\textwidth, clip=true]{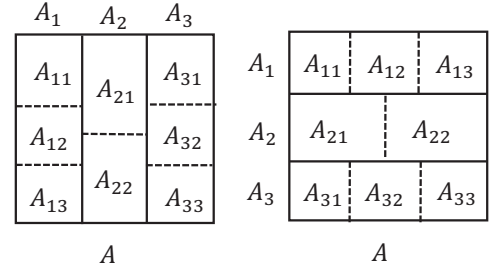}
\caption{Two heterogeneous partitions of matrix $A$: left graph illustrates Case 1 and right graph illustrates Case 2}
\end{figure}

\subsection{Heterogeneous partition of Case 1}
We first consider Case 1, where $A$ is partitioned into $c$ complete column blocks and each column block is partitioned into $c_{i}$ row subblocks, i.e.,
\begin{align}\label{b1}
A=[A_{1},A_{2}, \ldots, A_{c}], A_{i}=\text{col}(A_{i1}, A_{i2}, \ldots, A_{ic_{i}}),
\end{align}
with $A_{i} \in \mathbb{R}^{m\times n_{i}}, i\in \bm c$ and $A_{ij} \in \mathbb{R}^{m_{ij}\times n_{i}},j\in \bm c_{i}$. It follows that
\begin{align}\label{b2}
\sum^{c}_{i=1}n_{i}=n, \sum^{c_{i}}_{j=1}m_{ij}=m, i\in \bm c.
\end{align}
Accordingly, the vector $b$ is partitioned by
\begin{align*}
b=\sum^{c}_{i=1}b_{i}, b_{i}=\text{col}(b_{i1}, \ldots, b_{ic_{i}}), b_{ij}\in \mathbb{R}^{m_{ij}}.
\end{align*}
In this case, $b_{i}$ is chosen arbitrarily such that $b=\sum^{c}_{i=1}b_{i}$, but $b_{ij}$ is derived from a partition of $b_{i}$.

\begin{figure}[!ht]
\centering
\includegraphics[width=0.2\textwidth, clip=true]{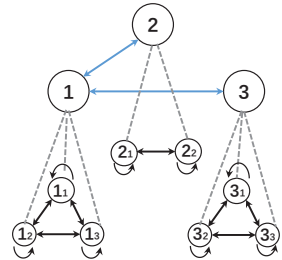}
\caption{Double-layered multi-agent network}
\end{figure}

Different from homogenous partition in \eqref{3} where single-layered grid multi-agent network is introduced, we consider a double-layered multi-agent network for heterogeneous partition, which is composed of $c$ clusters. Each cluster $i\in \bm c$ is composed of $c_{i}$ agents labelled by $i_{1},i_{2}, \ldots, i_{c_{i}}$. Each cluster does not involve with any computation, only collects and distributes all the agents information in the same cluster. The communication of all the clusters $\mathcal{V}=\{1,2,\ldots, c\}$ is described by the graph $\mathbb{G}=\{\mathcal{V},\mathcal{E}\}$ with the edge set $\mathcal{E}=\mathcal{V}\times \mathcal{V}$, and the neighbors of cluster $i$ is denoted by $\mathcal{N}_{i}=\{j\in \mathcal{V} | (j,i) \in \mathcal{E}\}$. The communication among all the agents in the same cluster $i$ is described by the graph $\mathbb{G}_{i}=\{\mathcal{V}_{i}, \mathcal{E}_{i}\}, i\in \bm c$, where $\mathcal{V}_{i}=\{i_{1},\ldots, i_{c_{i}}\}$ is the node set and $\mathcal{E}_{i}=\mathcal{V}_{i}\times \mathcal{V}_{i}$ is the edge set, and the neighbors of agent $i_{j}$ in the cluster $i$ is denoted by $\mathcal{N}_{ij}=\{i_{k}\in \mathcal{V}_{i} |  (i_{k},i_{j}) \in \mathcal{E}_{i}\}$. Fig. 4 shows such a double-layered multi-agent network, which corresponds to the example in left graph in Fig. 3.

\begin{Assumption}
The graph $\mathbb{G}$ is undirected and connected, and the graph $\mathbb{G}_{i}, i\in \bm c$ is undirected and connected.
\end{Assumption}

Suppose that each agent $i_{j}$ in the cluster $i$ knows $A_{ij} \in \mathbb{R}^{m_{ij}\times n_{i}}$ and $b_{ij} \in \mathbb{R}^{m_{ij}}$. In this section, we aim to develop a distributed algorithm such that each agent $i_{j}$'s state $x_{ij}\in \mathbb{R}^{n_{i}}$ exponentially converges to a constant vector $x^{*}_{ij}$ satisfying that

(i) All $x^{*}_{ij}$ in the same cluster $i$ achieve a consensus $x^{*}_{i}, i\in \bm c$, i.e.,
\begin{align}\label{b3}
x_{i1}=x_{i2}=\ldots=x_{ic_{i}}=x^{*}_{i}, i\in \bm c.
\end{align}

(ii) The consensus value $x^{*}_{i}$ of all the clusters satisfy that
\begin{align}\label{b4}
A^{T}_{i}\sum^{c}_{i=1}(A_{i}x^{*}_{i}-b_{i})=0, \forall i\in \bm c.
\end{align}
Note that $\sum^{c}_{i=1}(A_{i}x^{*}_{i}-b_{i})=Ax^{*}-b$ with $x^{*}=\text{col}(x^{*}_{1}, \ldots, $ $x^{*}_{c})$ and it follows from \eqref{b4} that $ A^{T}_{i}(Ax^{*}-b)=0,\forall i\in \bm c$, which implies that $A^{T}(Ax^{*}-b)=0$. Thus, we have that $x^{*}_{ij}$ collaboratively forms a least square solution of \eqref{1}.
\subsubsection{Distributed algorithm}
To ensure that all the agents' states $x_{ij}$ converge to $x^{*}_{ij}$ satisfying \eqref{b3} and \eqref{b4}, we introduce two extra variables $y_{ij}\in \mathbb{R}^{m_{ij}}, z_{ij}\in \mathbb{R}^{n_{i}}$ for each agent $i_{j}, i\in \bm c, j\in \bm c_{i}$. Then, a distributed algorithm for each agent $i_{j}$ is developed as
\begin{align}\label{b5}
\begin{cases}
\dot{x}_{ij}=-(A^{T}_{ij}y_{ij}-\sum_{i_{k}\in \mathcal{N}_{ij}}(z_{ij}-z_{ik}))\\
~~~~~~~~-\sum_{i_{k}\in \mathcal{N}_{ij}}(x_{ij}-x_{ik})\\
\dot{\xi}_{ij}=-\sum_{k\in \mathcal{N}_{i}}(y_{ij}-E_{ij}y_{k}), y_{ij}=\xi_{ij}\\
~~~~~~~~+A_{ij}x_{ij}-b_{ij}\\
\dot{z}_{ij}=A^{T}_{ij}y_{ij}-\sum_{i_{k}\in \mathcal{N}_{ij}}(z_{ij}-z_{ik})
\end{cases}
\end{align}
where $\xi_{ij} \in \mathbb{R}^{m_{ij}}$ is an intermediate variable and its initial value is selected as zero, $E_{ij}\in \mathbb{R}^{m_{ij}\times m}$ denotes a row partition of the identity matrix $I_{m}$ with $\text{col}(E_{i1}, \ldots, E_{ic_{i}})=I_{m}$. In particular, $E_{ij}\in \mathbb{R}^{m_{ij}\times m}$ is introduced such that $\frac{1}{c}\sum^{c}_{i=1}(\bar{A}_{i}\bar{x}_{i}-b_{i})$ can be estimated by the average consensus approach, where $\bar{A}_{i}=\text{diag}(A_{i1},\ldots, A_{ic_{i}}) \in \mathbb{R}^{m\times n_{i}c_{i}}$ and $\bar{x}_{i}=\text{col}(x_{i1},\ldots, x_{ic_{i}}) \in \mathbb{R}^{n_{i}c_{i}}$.

Define $\bar{y}_{i}=\text{col}(y_{i1}, \ldots, y_{ic_{i}})\in \mathbb{R}^{m}, \bar{z}_{i}=\text{col}(z_{i1},\ldots,$ $z_{ic_{i}})\in \mathbb{R}^{n_{i}c_{i}}$, $\bar{L}_{G_{i}}=L_{\mathbb{G}_{i}} \otimes I_{n_{i}}$ and $L_{\mathbb{\mathbb{G}}_{i}} \in \mathbb{R}^{c_{i}\times c_{i}}$ is the Laplacian matrix of the graph $\mathbb{G}_{i}$. Then, the algorithm \eqref{b5} is written as
\begin{align}\label{b6}
\begin{cases}
\dot{\bar{x}}_{i}=-\bar{A}^{T}_{i}\bar{y}_{i}+\bar{L}_{\mathbb{G}_{i}}(\bar{z}_{i}-\bar{x}_{i})\\
\dot{\bar{y}}_{i}=-\sum_{k\in \mathcal{N}_{i}}(\bar{y}_{i}-\bar{y}_{k})+\bar{A}_{i}\dot{\bar{x}}_{i}\\
\dot{\bar{z}}_{i}=\bar{A}^{T}_{i}\bar{y}_{i}-\bar{L}_{\mathbb{G}_{i}}\bar{z}_{i}
\end{cases}
\end{align}
Let $\hat{x}=\text{col}(\bar{x}_{1},\ldots, \bar{x}_{c})$, $\hat{y}=\text{col}(\bar{y}_{1},\ldots, \bar{y}_{c})$, $\hat{z}=\text{col}(\bar{z}_{1},$ $\ldots, \bar{z}_{c})$, $\hat{A}=\text{diag}(\bar{A}_{1},\ldots, \bar{A}_{c})$, $\hat{L}=\text{diag}(\bar{L}_{\mathbb{G}_{1}},\ldots, \bar{L}_{\mathbb{G}_{c}})$, $\hat{L}_{\mathbb{G}}=L_{\mathbb{G}}\otimes I_{n} $ and $L_{\mathbb{G}}\in \mathbb{R}^{c\times c}$ is the Laplacian matrix of the graph $\mathbb{G}$. The compact form of \eqref{b6} in terms of $(\hat{x},\hat{y},\hat{z})$ is given as
\begin{align}\label{b7}
\begin{cases}
\dot{\hat{x}}=-\hat{A}^{T}\hat{y}+\hat{L}\hat{z}-\hat{L}\hat{x},\\
\dot{\hat{y}}=-(\hat{L}_{\mathbb{G}}+\hat{A}\hat{A}^{T})\hat{y}+\hat{A}\hat{L}(\hat{z}-\hat{x}),\\
\dot{\hat{z}}=\hat{A}^{T}\hat{y}-\hat{L}\hat{z}.
\end{cases}
\end{align}
which is equivalent as
\begin{align}\label{b8}
\begin{bmatrix} \dot{\hat{x}} \\ \dot{\hat{y}} \\ \dot{\hat{z}}\end{bmatrix}=\begin{bmatrix} -\hat{L}_{\mathbb{G}}& -\hat{A}^{T}& \hat{L}_{\mathbb{G}}\\
-\hat{A}\hat{L}_{\mathbb{G}}& -(\hat{L}+\hat{A}\hat{A}^{T}) & \hat{A}\hat{L}_{\mathbb{G}} \\
0& \hat{A}^{T}& -\hat{L}_{\mathbb{G}} \end{bmatrix}\begin{bmatrix} \hat{x} \\ \hat{y} \\ \hat{z} \end{bmatrix}.
\end{align}

The following Theorem presents the convergence results of \eqref{b8}, and its proof is given in Appendix C.
\begin{Theorem}
Suppose that Assumption 2 holds and the initial value of $\xi_{ij}, i\in \bm c, j\in \bm c_{i} $ selected to be zero. The proposed distributed algorithm \eqref{b5} guarantees that all the agents's states $x_{ij}(t)$ converge exponentially to a constant value $x^{*}_{ij}$ satisfying \eqref{b3} and \eqref{b4}, which collaboratively forms a least square solution $x^{*}$ of equation \eqref{1}. Moreover, $x^{*}$ is an exact solution of \eqref{1} if and only if $\lim_{t\to \infty} \hat{y}(t)=0$.
\end{Theorem}

\subsection{Heterogeneous partition of Case 2}
In this section, we consider the heterogeneous partition of Case 2, where $A$ is partitioned into $c$ complete row blocks and each row block is partitioned into $c_{i}$ column subblocks, i.e.,
\begin{align}\label{a1}
A=\text{col}(A_{1},A_{2},\ldots,A_{c}), A_{i}=[A_{i1},A_{i2},\ldots,A_{ic_{i}}],
\end{align}
with $A_{i} \in \mathbb{R}^{m_{i}\times n}$ and $A_{ij} \in \mathbb{R}^{m_{i}\times n_{ij}}$. Then, it follows that
\begin{align}\label{a2}
\sum^{c}_{i=1}m_{i}=m, \sum^{c_{i}}_{j=1} n_{ij}=n, i\in \bm c.
\end{align}
Accordingly, the vector $b$ is partitioned by
\begin{align*}
b=\text{col}(b_{1}, b_{2}, \ldots, b_{c}), b_{i}=\sum^{c_{i}}_{j=1} b_{ij}, b_{ij}\in \mathbb{R}^{m_{i}},
\end{align*}
In this case, $b_{i}$ results from a partition of $b$ and $b_{ij}$ is chosen arbitrarily such that $b_{i}=\sum^{c_{i}}_{j=1} b_{ij}$.

Note that all the agents in the same cluster of Case 1 collectively know a complete column block. In this section, we also employ the double-layered multi-agent network $\mathbb{G}$ given in Fig. 4 while all the agents in the same cluster collectively know the complete row block. For the heterogenous partition of Case 2, suppose that each agent $i_{j}$ in the cluster $i$ knows $A_{ij} \in \mathbb{R}^{m_{i}\times n_{ij}}$ and $b_{ij} \in \mathbb{R}^{m_{i}}$. We aim to develop a distributed algorithm such that each agent $i_{j}$'s state $x_{ij} \in \mathbb{R}^{n_{ij}}$ converges exponentially to a constant vector $x^{*}_{ij}$ satisfying that

(i) All $\bar{x}^{*}_{i}=\text{col}(x^{*}_{i1}, x^{*}_{i2}, \ldots, x^{*}_{ic_{i}})\in \mathbb{R}^{n}$ among all the clusters reach a consensus $x^{*}$, i.e.,
\begin{align}\label{a3}
\bar{x}^{*}_{1}=\bar{x}^{*}_{2}=\ldots= \bar{x}^{*}_{c}=x^{*}.
\end{align}

(ii) The consensus value $x^{*}$ satisfies that
\begin{align}\label{a4}
\sum^{c}_{i=1}A^{T}_{i}(A_{i}x^{*}-b_{i})=0
\end{align}
Note from \eqref{a4} that $A^{T}(Ax^{*}-b)=0$, and we get that $x^{*}_{ij}$ collaboratively forms a least square solution of \eqref{1}.

\subsubsection{Distributed algorithm}
To ensure that all the agents' states $x_{ij}$ converge to $x^{*}_{ij}$ satisfying \eqref{a3} and \eqref{a4}, we develop a distributed algorithm for each agent $i_{j}$ is given as
\begin{align}\label{a5}
\begin{cases}
\dot{x}_{ij}=-(c_{i}A^{T}_{ij}y_{ij}-\sum_{k \in \mathcal{N}_{i}}(z_{ij}-E_{ij}z_{k}))\\
~~~~~~~~-\sum_{k \in \mathcal{N}_{i}}(x_{ij}-E_{ij}x_{k}),\\
\dot{\xi}_{ij}=-\sum_{i_{k} \in \mathcal{N}_{ij}}(y_{ij}-y_{ik}),y_{ij}=\xi_{ij}+\\
~~~~~~~~~(A_{ij}x_{ij}-b_{ij}),\\
\dot{z}_{ij}=c_{i}A^{T}_{ij}y_{ij}-\sum_{k \in \mathcal{N}_{i}}(z_{ij}-E_{ij}z_{k}),
\end{cases}
\end{align}
%
%
where the second equation of \eqref{a5} aims to estimate the average values $\frac{1}{c_{i}}\sum^{c_{i}}_{j=1}(A_{ij}x_{ij}(t)-b_{ij})$, and $E_{ij} \in \mathbb{R}^{n_{ij}\times n}$ denotes denotes a row partition of the identity matrix $I_{n}$ such that $\text{col}(E_{i1},E_{i2}, \ldots, E_{ic_{i}})=I_{n}$. Let $\bar{x}_{i}=\text{col}(x_{i1},\ldots, x_{ic_{i}}) \in \mathbb{R}^{n}, \bar{y}_{i}=\text{col}(y_{i1},\ldots, y_{ic_{i}})\in \mathbb{R}^{m_{i}c_{i}}, \bar{z}_{i}=\text{col}(z_{i1},\ldots, z_{ic_{i}})\in \mathbb{R}^{n}, \bar{A}_{i}=\text{diag}(A_{i1},\ldots, A_{ic_{i}})\in \mathbb{R}^{m_{i}c_{i}\times n}$, $\bar{L}_{G_{i}}=L_{\mathbb{G}_{i}} \otimes I_{m_{i}}\in \mathbb{R}^{m_{i}c_{i}\times m_{i}c_{i}}$ and $L_{\mathbb{\mathbb{G}}_{i}} \in \mathbb{R}^{c_{i}\times c_{i}}$ is the Laplacian matrix of the graph $\mathbb{G}_{i}$. Then, the algorithm \eqref{a5} is written as
\begin{align}\label{a6}
\begin{cases}
\dot{\bar{x}}_{i}=-c_{i}\bar{A}^{T}_{i}\bar{y}_{i}+\sum_{k\in \mathcal{N}_{i}}(\bar{z}_{i}-\bar{z}_{k}-(\bar{x}_{i}-\bar{x}_{k}))\\
\dot{\bar{y}}_{i}=-\bar{L}_{\mathbb{G}_{i}}\bar{y}_{i}+\bar{A}_{i}\dot{\bar{x}}_{i}\\
\dot{\bar{z}}_{i}=c_{i}\bar{A}^{T}_{i}\bar{y}_{i}-\sum_{k\in \mathcal{N}_{i}}(\bar{z}_{i}-\bar{z}_{k})
\end{cases}
\end{align}
Let  $\hat{x}=\text{col}(\bar{x}_{1},\ldots, \bar{x}_{c})$, $\hat{y}=\text{col}(\bar{y}_{1},\ldots, \bar{y}_{c})$, $\hat{z}=\text{col}(\bar{z}_{1},$ $\ldots, \bar{z}_{c})$, $\Gamma=\text{diag}(c_{1}I_{n}, \ldots, c_{c}I_{n})$,  $\hat{A}=\text{diag}(\bar{A}_{1},\ldots, \bar{A}_{c})$, $\hat{L}=\text{diag}(\bar{L}_{\mathbb{G}_{1}},\ldots, \bar{L}_{\mathbb{G}_{c}})$, $\hat{L}_{\mathbb{G}}=L_{\mathbb{G}}\otimes I_{n} $ and $L_{\mathbb{G}}\in \mathbb{R}^{c\times c}$ is the Laplacian matrix of the graph $\mathbb{G}$. A compact of \eqref{a6} in terms of $(\hat{x},\hat{y},\hat{z})$ is given as
\begin{align}\label{a7}
\begin{cases}
\dot{\hat{x}}=-\Gamma\hat{A}^{T}\hat{y}+\hat{L}_{\mathbb{G}}\hat{z}-\hat{L}_{\mathbb{G}}\hat{x},\\
\dot{\hat{y}}=-(\hat{L}+\hat{A}\Gamma\hat{A}^{T})\hat{y}+\hat{A}\hat{L}_{\mathbb{G}}(\hat{z}-\hat{x}),\\
\dot{\hat{z}}=\Gamma\hat{A}^{T}\hat{y}-\hat{L}_{\mathbb{G}}\hat{z}.
\end{cases}
\end{align}
which is equivalent as
\begin{align}\label{a8}
\begin{bmatrix} \dot{\hat{x}} \\ \dot{\hat{y}} \\ \dot{\hat{z}}\end{bmatrix}=\begin{bmatrix} -\hat{L}_{\mathbb{G}}& -\Gamma\hat{A}^{T}& \hat{L}_{\mathbb{G}}\\
-\hat{A}\hat{L}_{\mathbb{G}}& -(\hat{L}+\hat{A}\Gamma\hat{A}^{T}) & \hat{A}\hat{L}_{\mathbb{G}} \\
0& \Gamma\hat{A}& -\hat{L}_{\mathbb{G}} \end{bmatrix}\begin{bmatrix} \hat{x} \\ \hat{y} \\ \hat{z} \end{bmatrix}.
\end{align}
Since $\Gamma$ is added in the system \eqref{a8}, the results of Lemma 2 cannot be used directly. The following Lemma shows the eigenvalues of the system matrix of \eqref{a8}, and its proof is given in Appendix D.
\begin{Lemma}
All the non-zero eigenvalues of system matrix of \eqref{a8} have negative real part, and its zero eigenvalue must be non-defective.
\end{Lemma}

Based on the result of Lemma 3, the converge analysis of system \eqref{a8} is given in the following theorem, and its proof is given in Appendix E.
\begin{Theorem}
Suppose that Assumption 2 holds and the initial value of $\xi_{ij}$ selected to be zero. The proposed distributed algorithm \eqref{a5} guarantees that all the agents's states $x_{ij}(t)$ converge exponentially to a constant value $x^{*}_{ij}$ satisfying \eqref{a3} and \eqref{a4}, which collaboratively forms a least square solution $x^{*}$ of equation \eqref{1}. Moreover, $x^{*}$ is an exact solution of \eqref{1} if and only if $\lim_{t\to \infty} \hat{y}(t)=0$.
\end{Theorem}

\begin{Remark}
From the compact form \eqref{11}, \eqref{b8} and \eqref{a8} of the proposed three distributed algorithms \eqref{6}, \eqref{b5} and \eqref{a5} , we observe that the algorithm \eqref{6} has same system structure as \eqref{b5}, but the algorithm \eqref{a5} has different system structure with the others. The main reason is that the number of agents $c_{i}$ in the cluster $i$ is involved in \eqref{a5} but is not contained in the other algorithms. In fact, when $c_{i}$ is same for all $ i\in \bm c$, this parameter $c_{i}$ can be removed in \eqref{a5}, which also implies that the algorithm \eqref{6} can be regarded a special case of the algorithm \eqref{a5}.
\end{Remark}

\begin{Remark}
Different from consensus-based distributed algorithms in \cite{7,8,8a,8b} where all the agents' states have same dimension as the unknown variables, the proposed distributed algorithms \eqref{b5} and \eqref{a5} allow all the agents' states to have different dimensions that are much smaller than that of the unknown variables to be solved. This implies that the proposed algorithms \eqref{b5} and \eqref{a5} can be applied for heterogeneous agents with different storage and computation capability. In particular, when the number of the cluster and the number of the agents are sufficiently large, all the agents states are scalars, which is equivalent to the algorithm \eqref{6}.
 \end{Remark}

\subsection{Simulation}
In this section, we show the effectiveness of the proposed algorithms \eqref{b5} and \eqref{a5} by solving the linear equation $Ax=b$ in \eqref{s1}. Firstly, we consider the heterogeneous partition of Case 1 with $A_{11}=1,A_{12}=[2,1]^{T},A_{13}=2,A_{21}=[2,-1]^{T},A_{22}=[-2,2]^{T},A_{31}=[1,-1]^{T},A_{32}=4,A_{33}=-3$. In addition, $b_{1}=b, b_{2}=0_{4\times 1}, b_{3}=0_{4\times 1}$ is chosen such that $b=\sum^{3}_{i=1}b_{i}$, and $b_{ij}$ derives from a partition of $b_{i}$. We utilize the double-layered multi-agent network $\mathbb{G}$ in Fig. 4 to verify the distributed algorithm \eqref{b5}. Define $E_{2}(t)=\sum^{c}_{i=1}\sum^{c_{i}}_{j=1}\Vert x_{ij}(t)-x^{*}_{i}\Vert^{2}, x^{*}=\text{col}(x^{*}_{1},\ldots, x^{*}_{c})$ be the total estimated error and $Y_{e}(t)=\sum^{c}_{i=1}\sum^{c_{i}}_{j=1} \Vert y_{ij}(t)\Vert^{2}$. Fig. 5 shows that $E_{2}(t)$ converges exponentially fast to zero, which implies that all the states $x_{ij}(t)$ converge to the constant value $x^{*}_{i}$ as the $i$th sub-block of $x^{*}$. Fig. 5 also shows that $Y(t)$ converges a nonzero value, which implies that $\lim_{t\to \infty} \hat{y}(t) \ne 0$ and verifies that the solution $x^{*}$ is not an exact solution of \eqref{s1}.

\begin{figure}[!ht]
\centering
\includegraphics[width=0.45\textwidth, clip=true]{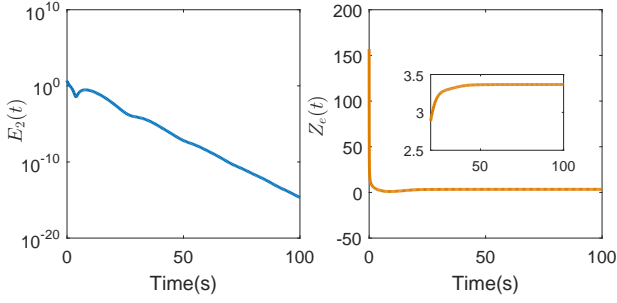}
\caption{Evolutions of $E_{2}(t)$ and $Y_{e}(t)$ under algorithm \eqref{b5}}
\end{figure}

We further consider the heterogeneous partition of Case 2 with $A_{11}=[1,2]^{T},A_{12}=[2,-1]^{T},A_{13}=[1,-1]^{T}, A_{21}=[1,-2],A_{22}=4, A_{31}=2, A_{32}=2, A_{33}=-2$.  In addition, $b_{i}, i=1,2,3$ is derived from a partition of $b$. $b_{11}=b_{1}, b_{22}=b_{2}, b_{33}=b_{3}$ are chosen and the others $b_{ij}$ are zeros such that $b_{i}=\sum^{c_{i}}_{j=1} b_{ij}$. The double-layered multi-agent network $\mathbb{G}$ in Fig. 4 is also utilized to demonstrate the distributed algorithm \eqref{a5}. Define $E_{3}(t)=\sum^{c}_{i=1}\Vert \bar{x}_{i}(t)-x^{*}\Vert^{2}, \bar{x}_{i}=\text{col}(x_{i1},\ldots, x_{ic_{i}})$ and $Y_{e}(t)=\sum^{c}_{i=1}\sum^{c_{i}}_{j=1}\Vert y_{ij}(t)\Vert^{2}$. Fig. 6 shows that $E_{3}(t)$ converges exponentially fast to zero, which implies that all the states $x_{ij}(t)$ converge to the constant value $x^{*}_{ij}$, and consists of a solution $x^{*}=\text{col}(x^{*}_{i1}, \ldots, x^{*}_{ic_{i}})$. Fig. 6 also shows that $Y_{e}(t)$ finally converges a non-zero value, which verifies that $x^{*}$ is not exact solution of \eqref{1}.

\begin{figure}[!ht]
\centering
\includegraphics[width=0.45\textwidth, clip=true]{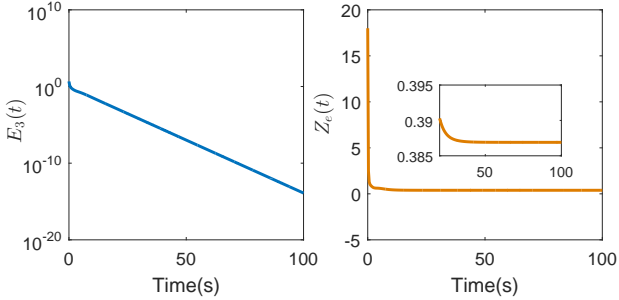}
\caption{Evolutions of $E_{1}(t)$ and $Y_{e}(t)$ under algorithm \eqref{a5}}
\end{figure}
%

\section{Conclusion}
This paper proposes three kinds of distributed algorithms via multi-agent networks to achieve a least square solution of linear equation under different partition structures. In particular, each agent has only access to a small block of linear equations, and only needs to control three scalar states when the linear equation is partitioned as small as possible. These advantages can be utilized to reduce the burdens of computation, storage and communication for solving the large-scale equations. Future works includes extending the proposed distributed algorithm to the time-varying network or a general directed graph.

\appendices
\section{Proof of Lemma 2}
Firstly, we prove the result of (i). Let $\lambda_{m} \in \mathbb{C}$ denote any eigenvalue of $M$, there exists a non-zero eigenvector $u=\text{col}(u_{1}, u_{2}, u_{3})$ in the complex vector space such that
\begin{align}\label{A1}
\begin{bmatrix} -L_{1} & -Y & L_{1}\\ -Y^{T}L_{1} & -Y^{T}Y-L_{2} & Y^{T}L_{1}\\ 0 & Y & -L_{1} \end{bmatrix}\begin{bmatrix} u_{1}\\u_{2}\\u_{3} \end{bmatrix}=\lambda_{m}\begin{bmatrix} u_{1}\\u_{2}\\u_{3} \end{bmatrix}.
\end{align}
We then have that
\begin{subequations}\label{A2}
\begin{align}
&-L_{1}u_{1}=\lambda_{m}(u_{1}+u_{3}), \label{A2.1}\\
&-L_{2}u_{2}=\lambda_{m}(-Y^{T}u_{1}+u_{2}), \label{A2.2}\\
&-L_{1}u_{3}=\lambda_{m}u_{3}-Yu_{2}. \label{A2.3}
\end{align}
\end{subequations}
Since $u$ belongs to the complex vector space, its conjugate transpose $u^{H}$ is introduced. Note that $Y^{T}=Y^{H}$ since $Y$ is a real matrix. By left multiplying $u^{H}_{2}$ of \eqref{A2.2}, one can derive that
\begin{align}\label{A3}
-u^{H}_{2}L_{2}u_{2}=-\lambda_{m}u^{H}_{2}Y^{H}u_{1}+\lambda_{m}u^{H}_{2}u_{2}.
\end{align}
Substituting \eqref{A2.3} and \eqref{A2.1} into \eqref{A3}, it follows that
\begin{align} \label{A4}
-u^{H}_{2}L_{2}u_{2}&=-\lambda_{m}(\lambda_{m}u_{3}+L_{1}u_{3})^{H}u_{1}+\lambda_{m}u^{H}_{2}u_{2} \notag \\
&=-\lambda_{m}u^{H}_{3}L_{1}u_{1}-\lambda_{m}\lambda^{H}_{m}u^{H}_{3}u_{1}+\lambda_{m}u^{H}_{2}u_{2} \notag \\
&=\lambda^{2}_{m}(u_{1}+u_{3})^{H}(u_{1}+u_{3})-\lambda^{2}_{m}u^{H}_{1}(u_{1}+u_{3}) \notag \\
&~~~-\lambda_{m}\lambda^{H}_{m}u^{H}_{3}u_{1}+\lambda_{m}u^{H}_{2}u_{2}.
\end{align}
In addition, by using \eqref{A2.1} again, one can derive that
\begin{subequations}\label{A5}
\begin{align}
&-\lambda^{2}_{m}u^{H}_{1}(u_{1}+u_{3})=\lambda_{m}u^{H}_{1}L_{1}u_{1},\\
&-\lambda_{m}\lambda^{H}_{m}u^{H}_{3}u_{1}=\lambda_{m}u^{H}_{1}L_{1}u_{1}+\lambda_{m}\lambda^{H}_{m}u^{H}_{1}u_{1}.
\end{align}
\end{subequations}
Substituting \eqref{A5} into \eqref{A4}, we have that
\begin{align}\label{A6}
\lambda^{2}_{m}(u_{1}&+u_{3})^{H}(u_{1}+u_{3})+\lambda_{m}(2u^{H}_{1}L_{1}u_{1}\notag \\
&+u^{H}_{2}u_{2})+\lambda_{m}\lambda^{H}_{m}u^{H}_{1}u_{1}+u^{H}_{2}L_{2}u_{2}=0.
\end{align}
Define $a=(u_{1}+u_{3})^{H}(u_{1}+u_{3}), b=2u^{H}_{1}L_{1}u_{1}+u^{H}_{2}u_{2}, c_{1}=u^{H}_{1}u_{1},c_{2}=u^{H}_{2}L_{2}u_{2}$. We easily know that $a,b,c_{1},c_{2}$ are all nonnegative real numbers. Then, equation \eqref{A6} is simplified as
\begin{align}\label{A7}
a\lambda^{2}_{m}+b\lambda_{m}+\lambda_{m}\lambda^{H}_{m}c_{1}+c_{2}=0.
\end{align}
Let $\lambda_{m}=\beta_{1}+\beta_{2}\bm i$, where $\beta_{1}, \beta_{2}$ are the real part and imaginary part, and $\bm i$ is imaginary unit. Then, \eqref{A7} can be written as
\begin{subequations}\label{A8}
\begin{align}
&a(\beta^{2}_{1}-\beta^{2}_{2})+\beta_{1}b+(\beta^{2}_{1}+\beta^{2}_{2})c_{1}+c_{2}=0, \label{A8.1}\\
&2a\beta_{1}\beta_{2}+\beta_{2}b=0 \label{A8.2}.
\end{align}
\end{subequations}
For any non-zero eigenvalue $\lambda_{m} \ne0$, the results of (i) is equivalent to show $\beta_{1}<0$. To this end, we consider two cases: $\beta_{2}=0$ for Case 1 and $\beta_{2}\ne 0$ for Case 2.

\textbf{Case 1:} by considering $\beta_{2}=0$, \eqref{A8} is arranged as
\begin{align}\label{A9}
(a+c_{1})\beta^{2}_{1}+\beta_{1}b+c_{2}=0.
\end{align}

(1) Firstly, we consider that $a+c_{1}=0$. It is obvious that $a=c_{1}=0$ since $a \ge 0,c_{1} \ge 0$. According to the expressions of $a, c$, we have that $u_{1}=0, u_{3}=0$. This implies that $u_{2}\ne 0$ since $u \ne0$, which further derive $b>0$. Thus, it follows that $\beta_{1}=-\frac{c_{2}}{b}\le 0$. Since $\lambda_{m}\ne 0$, we obtain that $\beta_{1}<0$.

(2) Secondly, we consider that $a+c_{1}\ne 0$. We first show $b>0$ by a a contradiction argument. Note that if $b=0$, we have $u_{2}=0$ and $c_{2}=0$. It follows from \eqref{A9} that $\beta_{1}=0$, which induces a contradiction with $\lambda_{m}=0$. Then, we have that $a+c>0,b>0$ and $c_{2}\ge 0$. According to the quadratic function root formula, it follows that the nonzero real root of \eqref{A9} is always negative when $a+c>0,b>0$ and $c_{2}\ge 0$. This implies that $\beta_{1}<0$.

\textbf{Case 2:} By considering $\beta_{2}\ne 0$, it follows from \eqref{A8.2} that
\begin{align}\label{A10}
2a\beta_{1}+b=0.
\end{align}
We first show that $a \ne 0$ by a contradiction argument. If $a=0$, it follows from \eqref{A10} that $b=0$. According to the expressions of $a,b$, we have that $u_{1}+u_{3}=0$ and $u_{2}=0$. In addition, substituting $a=0,b=0, u_{2}=0$ into \eqref{A8.1}, it follows that $(\beta^{2}_{1}+\beta^{2}_{2})c_{1}=0$. Since $\beta^{2}_{1}+\beta^{2}_{2}\ne 0$, one can derive that $c_{1}=0$ and $u_{1}=0$. This implies $u_{i}=0,i=1,2,3$, which induces a contradiction with $u\ne 0$. Then, we get that $a>0$. It follows from \eqref{A10} that $\beta_{1}=-\frac{b}{2a}\le 0$. Next, we show that $b>0$ by a contradiction argument. If $b=0$, we have that $L_{1}u_{1}=0$. It follows from $\eqref{A2.1}$ that $u_{1}+u_{3}=0$ since $\lambda_{m}\ne0$, which induces a contradiction with $a=(u_{1}+u_{3})^{H}(u_{1}+u_{3})\ne 0$. Since $b>0$, we obtain that $\beta_{1}<0$.

Secondly, we prove the result of (ii). By contradiction, we suppose that zero eigenvalue $\lambda_{m}=0$ of $M$ is defective, that is, its algebraic multiplicity is strictly greater than its geometric multiplicity. Then there exists a nonzero generalized eigenvector $v=\text{col}(v_{1},v_{2},v_{3})$ such that
\begin{align}\label{A11}
Mv=u, Mu=0.
\end{align}
Note that $Mu=0$ corresponds to $\lambda_{m}=0$ in \eqref{A1} and \eqref{A2}, and then we get that
\begin{align}\label{A12}
L_{1}u_{1}=0, L_{2}u_{2}=0, L_{1}u_{3}=Yu_{2}.
\end{align}
In addition, it follows from $Mv=u$ that
\begin{subequations}\label{A13}
\begin{align}
&-L_{1}v_{1}-Yv_{2}+L_{1}v_{3}=u_{1},\label{A13.1}\\
&-Y^{T}L_{1}v_{1}-(Y^{T}Y+L_{2})v_{2}+Y^{T}L_{1}v_{3}=u_{2},\label{A13.2}\\
&Yv_{2}-L_{1}v_{3}=u_{3}.\label{A13.3}
\end{align}
\end{subequations}
From \eqref{A13.1} and \eqref{A13.3}, it follows that $-L_{1}v_{1}=u_{1}+u_{3}$. By substituting $-L_{1}v_{1}=u_{1}+u_{3}$ and \eqref{A13.3} into \eqref{A13.2} yields
\begin{align}\label{A14}
Y^{T}u_{1}-L_{2}v_{2}=u_{2}.
\end{align}
From $Yu_{2}=L_{1}u_{3}$ and $L_{1}u_{1}=0$ in \eqref{A12}, we have that $u^{H}_{1}Yu_{2}=u^{H}_{1}L_{1}u_{3}=0$. By left multiplying $u^{H}_{2}$ of \eqref{A14}, one can derive that $u^{H}_{2}u_{2}=u^{H}_{2}Y^{H}u_{1}=0$. This implies that $u_{2}=0$, and then we derive that $L_{1}u_{3}=Yu_{2}=0$. Based on the facts that $-L_{1}v_{1}=u_{1}+u_{3}$, $L_{1}u_{3}=0$ and $L_{1}u_{1}=0$, we get that $(u_{1}+u_{3})^{H}(u_{1}+u_{3})=0$. This implies that $u_{1}+u_{3}=0$. In addition, by left multiplying $u^{H}_{3}$ of \eqref{A13.3}, it follows that $u^{H}_{3}Yv_{2}=u^{H}_{3}u_{3}$. Since $u_{1}=-u_{3}$ and $u_{2}=0$, \eqref{A14} is written as $Y^{T}u_{3}+L_{2}v_{2}=0$. By left multiplying $v^{H}_{2}$, we have that
\begin{align}\label{27}
v^{H}_{2}Y^{H}u_{3}+v^{H}_{2}L_{2}v_{2}=u^{H}_{3}u_{3}+v^{H}_{2}L_{2}v_{2}=0.
\end{align}
This implies that $u_{3}=0$ and $u_{1}=0$, which induces a contradiction with $u\ne 0$. Thus, we can conclude that zero eigenvalue of $M$ is non-defective.

\section{Proof of Theorem 1}
Since homogeneous partition \eqref{3} is regarded as a special case of heterogeneous partition of Case 1 in Section IV when $c=n$ and $c_{i}=m$ for any $i\in \bm c$. Thus, the proof of Theorem 1 is similar to that of Theorem 2, which is given in Appendix C.

\section{Proof of Theorem 2}
Firstly, we show that an equilibrium of system \eqref{b8} exists and is represented as $\tilde{\xi}^{*}=\text{col}(\tilde{x}^{*}, \tilde{y}^{*}, \tilde{z}^{*})$. Based on Lemma 1, we have that there must exists a constant vector $\eta \in \mathbb{R}^{n}$ such that $A^{T}(A\eta-b)=0$.  By virtue of heterogeneous partition of Case 1 in \eqref{b1} and $\eta=\text{col}(\eta_{1},\ldots, \eta_{c})$ with $\eta_{i} \in \mathbb{R}^{n_{i}}$, one can further derive that
\begin{align}\label{C1}
A^{T}_{i}\sum^{c}_{i=1}(A_{i}\eta_{i}-b_{i})=0.
\end{align}
Let $\tilde{y}^{*}=\text{col}(\tilde{y}^{*}_{1},\ldots, \tilde{y}^{*}_{c})=1_{c}\otimes \frac{1}{c}(A\eta-b)$ and $\tilde{x}^{*}=\text{col}(\tilde{x}^{*}_{1}, \ldots, \tilde{x}^{*}_{c})$ with $\tilde{x}^{*}_{i}=1_{c_{i}}\otimes \eta_{i}, i\in \bm c$. Under Assumption 2, we have that $\bar{L}_{\mathbb{G}_{i}}\tilde{x}^{*}_{i}=0, i\in \bm c$. It then follows that
\begin{align}\label{C2}
\hat{L}\tilde{x}^{*}=0, \hat{L}_{\mathbb{G}}\tilde{y}^{*}=0.
\end{align}
From the chosen $\tilde{y}^{*}$, we have that $\tilde{y}^{*}_{i}=\frac{1}{c}(A\eta-b)=\frac{1}{c}\sum^{c}_{i=1}(A_{i}\eta_{i}-b_{i})$, and it follows from \eqref{C1} that $A_{i}\tilde{y}^{*}_{i}=0$. By virtue of partition in \eqref{b1} and $\tilde{y}^{*}_{i}=\text{col}(\tilde{y}^{*}_{i1},\ldots, \tilde{y}^{*}_{ic_{i}})$ with $\tilde{y}^{*}_{ij} \in \mathbb{R}^{m_{ij}}$, we get that $\sum^{c_{i}}_{j=1}A_{ij}\tilde{y}_{ij}=0$, which can be written as
\begin{align}\label{C3}
(1^{T}_{c_{i}}\otimes I_{m_{i}}) \bar{A}^{T}_{i}\tilde{y}^{*}_{i}=0, \forall i\in \bm c.
\end{align}
Under Assumption 2, one has that $\text{Image}(\bar{L}_{\mathbb{G}_{i}})=\text{ker}(1^{T}_{c_{i}}\otimes I_{m_{i}})$. Then, there exists a constant vector $\tilde{z}^{*}_{i} \in \mathbb{R}^{m}$ such that
\begin{align}\label{C4}
\bar{A}^{T}_{i}\tilde{y}^{*}_{i}-\bar{L}_{\mathbb{G}_{i}}\tilde{z}^{*}_{i}=0, \forall i\in \bm c.
\end{align}
Define $\tilde{z}^{*}=\text{col}(\tilde{z}^{*}_{1},\ldots, \tilde{z}^{*}_{c})$, and it follows that
\begin{align}\label{C5}
\hat{A}^{T}\tilde{y}^{*}-\hat{L}\tilde{z}^{*}=0.
\end{align}
Based on \eqref{C2} and \eqref{C5}, we have that $\xi=\text{col}(\tilde{x}^{*},\tilde{y}^{*},\tilde{z})$ is an equilibrium of system \eqref{b8}.

Secondly, we analyze the convergence of the system state $(\hat{x}(t),\hat{y}(t),\hat{z}(t))$ with respect to the equilibrium $\tilde{\xi}^{*}$. Let $\xi(t)=\text{col}(\hat{x}(t),\hat{y}(t),\hat{z}(t))$ and the error $e(t)=\xi(t)-\tilde{\xi}^{*}$. Denote the system matrix of \eqref{b8} as $Q$, and one has that
\begin{align}\label{C6}
\dot{e}(t)=Qe(t).
\end{align}
From Lemma 2, we have that all the nonzero eigenvalues of $Q$ have negative real parts and zero eigenvalue is nondefective. Thus, there exists a constant vector $q \in \text{ker}(Q)$ such that $e(t)$ converges to $q$ with exponential convergence rate \cite{17}. Also, it follows that $\xi(t)$ converges exponentially to a constant vector $\hat{\xi}^{*}$, which satisfies that
\begin{align}\label{C7}
\hat{\xi}^{*}=\tilde{\xi}^{*}+q,
\end{align}
with $\hat{\xi}^{*}=\text{col}(\hat{x}^{*},\hat{y}^{*},\hat{z}^{*})$.

Finally, we only need to show that $\hat{x}^{*}=\text{col}(\bar{x}^{*}_{1}, \ldots, \bar{x}^{*}_{c})$ and $\bar{x}^{*}_{i}=\text{col}(x^{*}_{i1}, \ldots, x^{*}_{ic_{i}})$ satisfy \eqref{b3} and \eqref{b4}. From \eqref{C7}, one can derive that $\hat{\xi}^{*}$ is also an equilibrium of \eqref{b8}. Based on \eqref{b6}, we get that
\begin{subequations}\label{C8}
\begin{align}
\bar{A}^{T}_{i}\bar{y}^{*}_{i}-\bar{L}_{\mathbb{G}_{i}}\bar{z}^{*}_{i}+\bar{L}_{\mathbb{G}_{i}}\bar{x}^{*}_{i}=0, \label{C8.1}\\
 \sum_{k\in \mathcal{N}_{i}}(\bar{y}^{*}_{i}-\bar{y}^{*}_{k})=0,\label{C8.2}\\
 \bar{A}^{T}_{i}\bar{y}^{*}_{i}-\bar{L}_{\mathbb{G}_{i}}\bar{z}^{*}_{i}=0, \label{C8.3}
\end{align}
\end{subequations}
for any $i\in \bm c$. Under Assumption 2, from \eqref{C8.1} and \eqref{C8.3}, one has that $\bar{x}^{*}_{i}=1_{c_{i}} \otimes x^{*}_{i}$ with a constant vector $x^{*}_{i}$, which implies that \eqref{b3} is satisfied. In addition, from \eqref{C8.2}, it follows that $\bar{y}^{*}_{1}=\ldots=\bar{y}^{*}_{c}$. Moreover, with the initial condition $\xi_{ij}(0)=0, i\in \bm c, j\in \bm c_{i}$, one has that $\bar{y}_{i}(0)=\bar{A}_{i}\bar{x}_{i}(0)-b_{i}$. Summing the second equation \eqref{b6} from $i=1$ to $c$ can yield $\sum^{c}_{i=1}\dot{\bar{y}}_{i}(t)=\sum^{c}_{i=1}\bar{A}_{i}\dot{\bar{x}}_{i}$. Integrating it from $t_{0}$ to $t$, under the condition $\bar{y}_{i}(0)=\bar{A}_{i}\bar{x}_{i}(0)-b_{i}$, we have that $\sum^{c}_{i=1}\bar{y}_{i}(t)=\sum^{c}_{i=1}(\bar{A}_{i}\bar{x}_{i}(t)-b_{i})$ always holds for $\forall t\ge 0$. Noting that $\bar{y}^{*}_{1}=\ldots=\bar{y}^{*}_{c}$ and $\bar{x}^{*}_{i}=1_{c_{i}} \otimes x^{*}_{i}$, one can derive that $\bar{y}^{*}_{i}=\frac{1}{c}\sum^{c}_{i=1}(A_{i}x^{*}_{i}-b_{i}), \forall i\in \bm c$. By left multiplying $1^{T}_{c_{i}} \otimes I_{n_{i}}$ of \eqref{C8.1}, we have that $(1^{T}_{c_{i}} \otimes I_{n_{i}})\bar{A}^{T}_{i}\bar{y}^{*}_{i}=A^{T}_{i}\bar{y}^{*}_{i}=\frac{1}{c}A^{T}_{i}\sum^{c}_{i=1}(A_{i}x^{*}_{i}-b_{i})=\frac{1}{c}A^{T}_{i}(Ax^{*}-b)=0$. This implies that \eqref{b4} holds. Thus, we have shown that $x_{ij}(t)$ converge exponentially to $x^{*}_{ij}$, which collaboratively forms a least square solution of equation \eqref{1}.

In addition, from above results, we have obtained that $\hat{y}(t)$ converges exponentially to $\hat{y}^{*}=1_{c}\otimes \frac{1}{c}(Ax^{*}-b)$. It is clear that $x^{*}$ is an exact solution of \eqref{1} if and only if $\hat{y}^{*}=0$. It is indicated that the final convergence values of $y_{ij}(t)$ can be employed to verify whether the obtained solution $x^{*}$ is an exact solution.

\section{Proof of Lemma 3}
Denote the system matrix of system \eqref{a8} as $Q$, and we first show that all the non-zero eigenvalues $Q$ have negative real part. With the same eigenvalue and eigenvector notion in the proof of Lemma 2, it follows from $Q$ that
\begin{subequations}\label{B1}
\begin{align}
&-\hat{L}_{\mathbb{G}}u_{1}=\lambda_{m}(u_{1}+u_{3}), \label{B1.1}\\
&-\hat{L}u_{2}=\lambda_{m}(-\hat{A}u_{1}+u_{2}),\label{B1.2}\\
&T\hat{A}^{T}u_{2}-\hat{L}_{\mathbb{G}}u_{3}=\lambda_{m}u_{3}. \label{B1.3}
\end{align}
\end{subequations}
Let $u_{1}=\text{col}(\bar{u}_{1i},\ldots, \bar{u}_{1c}), u_{2}=\text{col}(\bar{u}_{2i},\ldots, \bar{u}_{2c})$ and $u_{3}=\text{col}(\bar{u}_{3i},\ldots, \bar{u}_{3c})$. From \eqref{B1.2} and \eqref{B1.3}, it follows that
\begin{subequations}\label{B2}
\begin{align}
&-\bar{L}_{\mathbb{G}_{i}}\bar{u}_{2i}=\lambda_{m}(-\bar{A}_{i}\bar{u}_{1i}+\bar{u}_{2i}), \label{B2.1}\\
&c_{i}\bar{A}^{T}_{i}\bar{u}_{2i}-\sum_{j\in \mathcal{N}_{i}}(\bar{u}_{3i}-\bar{u}_{3j})=\lambda_{m}\bar{u}_{3i}, \label{B2.2}
\end{align}
\end{subequations}
with $i\in \bm c$. By left multiplying $\bar{u}^{T}_{2i}$ of \eqref{B2.1}, substituting \eqref{B2.2} into \eqref{B2.1} can yield
\begin{align*}
&-c_{i}\bar{u}^{H}_{2i}\bar{L}_{\mathbb{G}_{i}}\bar{u}_{2i}=-\lambda_{m}c_{i}\bar{u}^{H}_{2i}\bar{A}_{i}\bar{u}_{1i}+\lambda_{m}c_{i}\bar{u}^{H}_{2i}\bar{u}_{2i}\\
&=-\lambda_{m}(\sum_{j\in \mathcal{N}_{i}}(\bar{u}_{3i}-\bar{u}_{3j}))^{H}\bar{u}_{1i}-\lambda_{m}\lambda^{H}_{m}\bar{u}^{H}_{3i}\bar{u}_{1i}+\lambda_{m}c_{i}\bar{u}^{H}_{2i}\bar{u}_{2i}.
\end{align*}
By summing above equation from $i=1$ to $c$, it follows that
\begin{align*}
-u^{H}_{2}\tilde{L}u_{2}=-\lambda_{m}u^{H}_{3}\hat{L}_{\mathbb{G}}u_{1}-\lambda_{m}\lambda^{H}_{m}u^{H}_{3}u_{1}+\lambda_{m}u^{T}_{2}\Xi u_{2},
\end{align*}
with $\tilde{L}=\text{diag}(c_{1}\bar{L}_{\mathbb{G}_{1}},\ldots, c_{c}\bar{L}_{\mathbb{G}_{c}})\ge 0$ and $\Xi=\text{diag}(c_{1}I_{c_{1}m_{1}},$ $\ldots, c_{c}I_{c_{c}m_{c}})>0$. By using \eqref{B1.1} and \eqref{B2.2}, we have that
\begin{align}\label{B3}
\lambda^{2}_{m}(u_{1}&+u_{3})^{H}(u_{1}+u_{3})+\lambda_{m}(2u^{H}_{1}\hat{L}_{\mathbb{G}}u_{1}+u^{H}_{2}\Xi u_{2}) \notag \\
&+\lambda_{m}\lambda^{H}_{m}u^{H}_{1}u_{1}+u^{H}_{2}\tilde{L}u_{2}=0.
\end{align}
Note that \eqref{B3} is almost same as \eqref{A6}. Thus, similar to the analyse of the Lemma 2, we easily get that
all the nonzero eigenvalues of $Q$ have negative real part.

Next, we show that $\lambda_{m}=0$ of $Q$ is defective. The proof process is similar to that of Lemma 2 and here only some crucial steps are given. Similar to Lemma 2, we have that $Mv=u$ and $Mu=0$. It follows from $Mu=0$ that $\hat{L}_{G}u_{1}=0, \hat{L}u_{2}=0, T\hat{A}^{T}u_{2}=\hat{L}_{\mathbb{G}}u_{3}$. In addition, from $Mv=u$, we have that
\begin{subequations}\label{B4}
\begin{align}
&-\hat{L}_{\mathbb{G}}v_{1}-T\hat{A}^{T}v_{2}+\hat{L}_{\mathbb{G}}v_{3}=u_{1},\label{B4.1}\\
&-\hat{A}\hat{L}_{\mathbb{G}}v_{1}-(\hat{A}T\hat{A}^{T}+\hat{L})v_{2}+\hat{A}\hat{L}_{\mathbb{G}}v_{3}=u_{2},\label{B4.2}\\
&T\hat{A}^{T}v_{2}-\hat{L}_{\mathbb{G}}v_{3}=u_{3}.\label{B4.3}
\end{align}
\end{subequations}
From above \eqref{B4}, we get that $-\hat{L}_{\mathbb{G}}v_{1}=u_{1}+u_{3}$ and $\hat{A}u_{1}-\hat{L}_{2}v_{2}=u_{2}$. This implies that  $c_{i}\bar{A}_{i}\bar{u}_{1i}-c_{i}\bar{L}_{\mathbb{G}_{i}}\bar{v}_{2i}=c_{i}\bar{u}_{2i}, \forall i\in \bm c$, which is arranged as
\begin{align}\label{B5}
\hat{A}Tu_{1}-\tilde{L}v_{2}=\Xi u_{2}.
\end{align}
 Since $\hat{L}u_{2}=0, T\hat{A}^{T}u_{2}=\hat{L}_{\mathbb{G}}u_{1}$ and $T^{H}=T$, we have that $\tilde{L}u_{2}=0$ and $u^{H}_{1}T\hat{A}^{T}u_{2}=0$. By left multiplying $u^{H}_{2}$ of \eqref{B5}, we have that $u^{H}_{2}\Xi u_{2}=0$, which implies that $u_{2}=0$. Then, we also get that
\begin{align}\label{B6}
u_{1}=-u_{3}, u^{H}_{3}u_{3}+v^{H}_{2}\tilde{L}v_{2}=0.
\end{align}
This implies that $u_{1}=u_{3}=0$, which occurs a contradiction. Then we conclude that $\lambda_{m}=0$ is defective.

\section{Proof of Theorem 3}
The proof process of Theorem 3 is similar to that of Theorem 2, and only some crucial steps are provided. Firstly, with the same notation in Theorem 2, we show that there exists a constant vector $\tilde{\xi}^{*}=\text{col}(\tilde{x}^{*}, \tilde{y}^{*}, \tilde{z}^{*})$ as an equilibrium of system \eqref{a8}. Since there exists $\eta \in \mathbb{R}^{n}$ satisfying $A^{T}(A\eta-b)=0$, we choose that $\tilde{x}^{*}=1^{T}_{c}\otimes \eta, \tilde{y}=\text{col}(\tilde{y}^{*}_{1},\ldots, \tilde{y}^{*}_{c})$ with $\tilde{y}^{*}_{i}=1^{T}_{c_{i}} \otimes \frac{1}{c_{i}}(A_{i}\eta-b_{i}), i \in \bm c$. Under Assumption 2, it follows that
\begin{align}\label{D1}
\hat{L}_{\mathbb{G}}\tilde{x}^{*}=0, \hat{L}\tilde{y}^{*}=0.
\end{align}
According to heterogeneous partition in \eqref{a1} and the chosen $\tilde{y}^{*}$, we have that $A^{T}(A\eta-b)=\sum^{c}_{i=1}A^{T}_{i}(A_{i}\eta-b_{i})=\sum^{c}_{i=1}c_{i}\bar{A}^{T}_{i}\tilde{y}^{*}_{i}=0$, and the last equation is written as
\begin{align}\label{D2}
(1^{T}_{c} \otimes I_{n})\Gamma \hat{A}^{T}\tilde{y}^{*}=0.
\end{align}
Since $\text{ker}(1^{T}\otimes I_{n})=\text{Image}(\hat{L}_{\mathbb{G}})$, then there exists a constant vector $\tilde{z}^{*}\in \mathbb{R}^{nc\times nc}$ such that
\begin{align}\label{D3}
\Gamma\hat{A}^{T}\tilde{y}^{*}-\hat{L}_{\mathbb{G}}\tilde{z}^{*}=0.
\end{align}
From \eqref{D1} and \eqref{D1}, one can derive that $\tilde{\xi}^{*}$ is an equilibrium of system \eqref{a8}.

Similar to the proof of Theorem 2, we have that $\hat{\xi}(t)$ converges exponentially fast to a constant vector $\hat{\xi}^{*}=\text{col}(\hat{x}^{*},\hat{y}^{*},\hat{z}^{*})$ satisfying that $\hat{\xi}^{*}=\tilde{\xi}^{*}+q$ with $q \in \text{ker}(Q)$. Next, we show that $\hat{x}^{*}=\text{col}(\bar{x}^{*}_{1},\ldots, \bar{x}^{*}_{c})$ satisfies \eqref{a3} and \eqref{a4}. Note that $\hat{\xi}^{*}$ is also an equilibrium of \eqref{a8}, and it then follows that
\begin{subequations}\label{D4}
\begin{align}
\sum_{k\in \mathcal{N}_{i}}(\bar{x}^{*}_{i}-\bar{x}^{*}_{k})=0,\label{D4.1}\\
\bar{L}_{\mathbb{G}_{i}}\bar{y}^{*}_{i}=0,\label{D4.2}\\
c_{i}\bar{A}^{T}_{i}\bar{y}^{*}_{i}-\sum_{k\in \mathcal{N}_{i}}(\bar{z}^{*}_{i}-\bar{z}^{*}_{k})=0. \label{D4.3}
\end{align}
\end{subequations}
for $i\in \bm c$. Under Assumption 2, it follows from \eqref{D4.1} that $\bar{x}^{*}_{1}=\ldots=\bar{x}^{*}_{c}=x^{*}$ with a constant vector $x^{*}$, which implies that \eqref{a3} holds. In addition, from \eqref{D4.2}, we get that $y^{*}_{i1}=\ldots=y^{*}_{ic_{i}}$. Under the initial condition $\xi_{ij}(0)=0$, similar to the analysis in Theorem 2, we have that $\sum^{c_{i}}_{j=1}y_{ij}(t)=\sum^{c_{i}}_{j=1}(A_{ij}x_{ij}(t)-b_{ij})$ always holds for $\forall t\ge 0$, and one further derive that $y^{*}_{ij}=\frac{1}{c_{i}}\sum^{c_{i}}_{j=1}(A_{ij}x^{*}_{ij}-b_{i})=\frac{1}{c_{i}}(A_{i}x^{*}-b_{i})$ and $\bar{y}^{*}_{i}=1_{c_{i}} \otimes \frac{1}{c_{i}}(A_{i}x^{*}-b_{i})$. From \eqref{D4.3}, one has that $\sum^{c}_{i=1}c_{i}\bar{A}^{T}_{i}\bar{y}^{*}_{i}=\sum^{c}_{i=1}A^{T}_{i}(A_{i}x^{*}-b)=0$, which implies that \eqref{a4}.  Thus, we have shown that $x_{ij}(t)$ converge exponentially to $x^{*}_{ij}$. In addition, we obtain that $\hat{y}(t)$ converges exponentially to $\hat{y}^{*}=\text{col}(\bar{y}^{*}_{1}, \ldots, \bar{y}^{*}_{c})$ with $\bar{y}_{i}=1_{c_{i}} \otimes \frac{1}{c_{i}}(A_{i}x^{*}-b_{i}), i\in \bm c$. Then we have that $x^{*}$ is an exact solution if and only if $\hat{y}^{*}=0$.

\ifCLASSOPTIONcaptionsoff
  \newpage
\fi

\end{document}